# New Analytical Expressions for the Levi-Civita Symbol and Its Treatment as a Generalized Function

W. Astar
University of Maryland, Baltimore County(UMBC); Baltimore, Maryland 21250

***Abstract -*** New, analytical expressions are found for the Levi-Civita symbol using the Kronecker delta symbol. The expressions are derived up to 3 dimensions, extended to higher dimensions, and confirmed in Matlab for 5 dimensions. The expressions can be re-cast in terms of elementary and/or special functions, which lead to the conclusion that the Levi-Civita Symbol can be treated as a generalized, discrete function

## 1. Introduction

The Levi-Civita symbol (or epsilon) is ascribed to Tullio Levi-Civita (1873-1941), who was an Italian physicist and a mathematician, and a contemporary of Albert Einstein. The epsilon is widely used in linear algebra, tensor analysis, and differential geometry. It represents the set {0, 1, -1}, and is generally defined for an *N*-dimensional vector space ($R^N$), as the sign of a permutation of the set of natural numbers {1, 2, 3, $\cdots$, *N*}. Using the index set or *N*-tuple {$i_1$, $i_2$, $i_3$, $\cdots$, $i_N$}, with each index $i_n$ assuming the values of the set {1, 2, 3,$\cdots$, *N*}, the epsilon is generally expressed as [1]:

$$\varepsilon_{i_1 i_2 i_3 \cdots i_N} = \begin{cases} +1 & \text{if } \{i_1,i_2,i_3,\cdots,i_N\} \text{ is any even permutation of } \{1,2,3,\cdots,N\}, \\ -1 & \text{if } \{i_1,i_2,i_3,\cdots,i_N\} \text{ is any odd permutation of } \{1,2,3,\cdots,N\}, \\ 0 & \text{otherwise.} \end{cases} \quad (1.1)$$

In $R^N$, an epsilon generates $N^N$ distinct *N*-tuples, out of which (*N*!) *N*-tuples have non-repetitive numbers, and ($N^N$-*N*!) *N*-tuples have at least one repetitive number. Expressed differently, the epsilon yields +1 for a set of (*N*!/2) *N*-tuples, -1 for another, different set of (*N*!/2) *N*-tuples, and 0 for the remaining ($N^N$-*N*!) *N*-tuples. Interchanging any 2 numbers in an *N*-tuple results in a sign change for the epsilon.

In two-dimensions ($R^2$), the epsilon is defined as

$$\varepsilon_{ij} = \begin{cases} +1 & \text{if } \{i,j\} = \{1,2\}, \\ -1 & \text{if } \{i,j\} = \{2,1\}, \\ 0 & \text{if } \{i,j\} \in \{\{1,1\},\{2,2\}\}, \end{cases} \quad (1.2)$$

whereas for three-dimensions ($R^3$),

$$\varepsilon_{ijk} = \begin{cases} +1 & \text{if } \{i,j,k\} \in \{\{1,2,3\},\{2,3,1\},\{3,1,2\}\}, \\ -1 & \text{if } \{i,j,k\} \in \{\{1,3,2\},\{2,1,3\},\{3,2,1\}\}, \\ 0 & \text{if } i=j, \text{ or } i=k, \text{ or } j=k. \end{cases} \quad (1.3)$$

It is difficult to view the above relations as conventional equations, in the same sense as well-known equations such as Newton's Second Law $F = m \cdot a$, or Einstein's mass-energy equivalence $E = m \cdot c^2$. Indeed, the above relations are more akin to look-up tables, or to a collection of conditions, than to an equation in which values can be substituted into variables. This difficulty is exacerbated for higher dimensions, and motivates the search for an alternative compact formula.

Some attempts have been made at presenting the epsilon as an equation. In $R^N$, a general expression is known to be [1]

$$\varepsilon_{i_1 i_2 i_3 \cdots i_N} = \prod_{N \geq n \geq m \geq 1} \mathrm{sgn}(i_n - i_m); \quad i_1, i_2, i_3, \cdots, i_N \in \{1, 2, 3, \cdots, N\} \tag{1.4}$$

with the signum (or the sign) function *generally* defined as

$$\mathrm{sgn}(x - y) = \frac{x - y}{|x - y|}; \quad x, y \in R. \tag{1.5}$$

The number of terms generated by the product (or pi) symbol in (1.4) is given by

$$\sum_{n=1}^{N-1} n = \frac{N(N-1)}{2}. \tag{1.6}$$

The relation (1.4) is relatively compact, and is easily generalized to any number of dimensions $N$. It is especially useful when the number of dimensions exceeds four. Although clearly a significant improvement over (1.1-3), it still incorporates an unresolved condition in the form of a multiple inequality, $N \geq n \geq m \geq 1$. Furthermore, since it has no anti-derivative due to a discontinuity at $x = y$, the signum function (1.5) is not an elementary function according to the Risch Algorithm [2-4].

In terms of other non-elementary functions, Straub [5] offered the following general expression for $R^N$, in terms of a $N \times N$ determinant of Kronecker delta symbols (or deltas),

$$\varepsilon_{i_1 i_2 i_3 \cdots i_N} = \begin{vmatrix} \delta_{i_1-1} & \delta_{i_1-2} & \delta_{i_1-3} & \cdots & \delta_{i_1-N} \\ \delta_{i_2-1} & \delta_{i_2-2} & \delta_{i_2-3} & \cdots & \delta_{i_2-N} \\ \delta_{i_3-1} & \delta_{i_3-2} & \delta_{i_3-3} & \cdots & \delta_{i_3-N} \\ \vdots & \vdots & \vdots & \ddots & \vdots \\ \delta_{i_N-1} & \delta_{i_N-2} & \delta_{i_N-3} & \cdots & \delta_{i_N-N} \end{vmatrix} \tag{1.7}$$

which can be specialized to any $R^N$. In $R^2$ for instance, it yields the simple relation [5],

$$\varepsilon_{ij} = \begin{vmatrix} \delta_{i-1} & \delta_{i-2} \\ \delta_{j-1} & \delta_{j-2} \end{vmatrix} = \delta_{i-1}\delta_{j-2} - \delta_{i-2}\delta_{j-1}. \tag{1.8}$$

In such applications, the delta symbol [6] due to L. Kronecker is variously defined as:

$$\delta_{z-z_0} = \delta_{zz_0} = \delta_{z_0}^z = \delta[z - z_0] = \begin{cases} 1, & z = z_0 \\ 0, & z \neq z_0 \end{cases} ; \; z_0 \in \mathbb{Z}. \tag{1.9}$$

Using the above definition, composite relations can be obtained, such as the well-known relation in $R^2$[1],

$$\varepsilon_{ij}\varepsilon_{mn} = \delta_{im}\delta_{jn} - \delta_{in}\delta_{jm} ; \; i, j, m, n \in \{1, 2\}. \tag{1.10}$$

It appears that the epsilon is not a heavily investigated concept, as Levi-Civita's original definition (1.1-3) had gained wide acceptance since its introduction ≈ 100 years ago. The wide availability of computers has also facilitated its implementation, especially critical when $N$ exceeds 3 or 4. Regardless, there have been some attempts at analytical alternatives to (1.1) and (1.4), especially for lower $N$. A recent instance is due to Jaramillo, who proffered without a derivation, in $R^3$, the expression [7]

$$\varepsilon_{ijk} = \frac{1}{2}(i-j)(j-k)(k-i) ; \; i, j, k \in \{1, 2, 3\} \tag{1.11}$$

which is certainly simpler, more succinct, and therefore more memorable than either (1.3), or (1.4) when specialized to $N = 3$. The expression is also in terms of elementary functions, unlike (1.4). However, Jaramillo made no attempt at generalizing (1.11) to higher dimensions.

In this report, analytical expressions are derived for the epsilon in $R^2$, and in $R^3$. Although deriving similar expressions for higher dimensional space is not within the scope of this report, they are deduced from those for the lower dimensional spaces, and verified in Matlab. In $R^N$ in general, it is found that

$$\varepsilon_{i_1 i_2 i_3 \cdots i_N} = \prod_{m=1}^{N-1} \prod_{n=1}^{N-m} \frac{i_{N+1-m} - i_n}{N+1-m-n} ; \; i_1, i_2, i_3, \cdots, i_N \in \{1, 2, 3, \cdots, N\}, \; N > 1, \tag{1.12}$$

is the simplest expression, and wholly in terms of elementary functions.

The epsilon as a generalized function is also explored in this report. It is found that many elementary and special functions $G$ may be utilized in a general equation (1.13), that can be used to replicate the results of an epsilon in any $R^N$ space,

$$\varepsilon_{i_1 i_2 i_3 \cdots i_N}(\lambda) = \prod_{m=1}^{N-1} \prod_{n=1}^{N-m} \frac{G(i_{N+1-m}\lambda) - G(i_n\lambda)}{G((N+1-m)\lambda) - G(n\lambda)} ; \; i_1, i_2, i_3, \cdots, i_N \in \{1, 2, 3, \cdots, N\}; N > 1; \lambda \in \mathbb{C}.$$

$$\tag{1.13}$$

## 2. Two dimensional vector space

In two-dimensions ($R^2$), the Levi-Civita symbol, or the epsilon, is defined as

$$\varepsilon_{ij} = \begin{cases} +1 & \text{if } \{i,j\} = \{1,2\}, \\ -1 & \text{if } \{i,j\} = \{2,1\}, \\ 0 & \text{if } \{i,j\} \in \{\{1,1\},\{2,2\}\}. \end{cases} \tag{2.1}$$

There are a total of $2^2 = 4$ possibilities, 2! of which yield non-zero results. This epsilon is also frequently expressed as

$$\varepsilon_{ij} = \frac{j-i}{|j-i|} = \text{sgn}(j-i); \quad i,j \in \{1,2\}, \tag{2.2}$$

which is a well-known expression [1], although its original provenance is obscure. This expression is amenable to a generalization to higher-order epsilons. Eq. (2.1) may also be configured in terms of a *signed* logical exclusive-OR,

$$\varepsilon_{ij} = (-1)^j \cdot (i-1) \oplus (j-1); \quad i,j \in \{1,2\}. \tag{2.3}$$

A procedure will now be outlined in this section that can be adapted to higher-order epsilons. Expressing the epsilon in terms of inter-indicial conditions, after arbitrarily declaring $i$ as the independent index with a domain of $\{1,2\}$,

$$\varepsilon_{ij} = \begin{cases} +1 & \text{if } j = (i+1)\delta_{i-1} = 2\delta_{i-1}, \\ -1 & \text{if } j = (i-1)\delta_{i-2} = \delta_{i-2}, \\ 0 & \text{if } j = i. \end{cases} \tag{2.4}$$

In terms of the Kronecker delta symbol (or the delta), the expression reduces to

$$\varepsilon_{ij} = \delta_{j-2\delta_{i-1}} - \delta_{j-\delta_{i-2}}. \tag{2.5}$$

However, since $i$ and $j$ are clearly constrained to the set $\{1, 2\}$, it is deduced that

$$\delta_{i-1} = 2 - i, \tag{2.6}$$

$$\delta_{i-2} = i - 1, \tag{2.7}$$

which collectively lead to the simplification that

$$\varepsilon_{ij} = \delta_{j-(4-2i)} - \delta_{j-(i-1)}. \tag{2.8}$$

This expression is simpler than (1.8) due to Straub [6]. According to (2.1), $i$ and $j$ are both constrained to $\{1, 2\}$. Thus, the 1st term is only non-zero when $i = 1$ and $j = 2$, or $(i - 2)(1 - j)$. Furthermore, the 2nd term is only non-zero when $i = 2$ and $j = 1$, and is thus equivalent to $(i - 1)(2 - j)$. Therefore

$$\varepsilon_{ij} = (i-2)(1-j) - (i-1)(2-j) = j - i \tag{2.9}$$

which is equivalent to (2.2), although simpler. This expression is also in terms of elementary functions, in accordance with the Risch Algorithm, unlike (2.2). It is clearly a trivial relation that should be obvious from a cursory inspection of (2.1). However, the approach taken to attain (2.9) will be adapted to higher dimensions, and is carried out here just to illustrate its utility for the simplest of epsilons, which is (2.1).

It appears that (2.9) is merely the simplest case of a more general relation. Many elementary and non-elementary functions $G$ satisfy (2.1), but in the form of

$$\varepsilon_{ij}(\lambda) = \frac{G(i\lambda) - G(j\lambda)}{G(\lambda) - G(2\lambda)}; \ i, j \in \{1, 2\} \ \& \ \lambda \in \mathbf{C}, \tag{2.10}$$

with $\mathbf{C}$ being the set of complex numbers. In terms of the identity function $G(\xi) = \xi$,

$$\varepsilon_{ij}(\lambda) = \frac{i\lambda - j\lambda}{\lambda - 2\lambda} = j - i, \tag{2.11}$$

which is (2.9), and without restrictions on $\lambda$. In terms of a trigonometric function, which is another elementary function,

$$\varepsilon_{ij}(\lambda) = \frac{\cos(i\lambda) - \cos(j\lambda)}{\cos(\lambda) - \cos(2\lambda)}, \ \lambda \in \mathbf{C} \tag{2.12}$$

with the simplest instance being that for $\lambda = \pi / 2$,

$$\varepsilon_{ij}(\pi/2) = \cos(i\pi/2) - \cos(j\pi/2) = \sin((j-i)\pi/2). \tag{2.13}$$

Clearly, other trigonometric functions are also viable in (2.10).

It is also feasible to re-express (2.9) using (2.10) with $G(\xi) = \Gamma(\xi)$,

$$\varepsilon_{ij}(\lambda) = \frac{\Gamma(i\lambda) - \Gamma(j\lambda)}{\Gamma(\lambda) - \Gamma(2\lambda)}, \ \lambda \in \mathbf{C}, \tag{2.14}$$

which is in terms of the gamma function, a non-elementary function. The constant $\lambda$ is still arbitrarily complex, as is the case for (2.11) and (2.12). However, the simplest case occurs for $\lambda = 1$ and for $G(\xi) = \Gamma(\xi + 1)$,

$$\varepsilon_{ij}(1) = \Gamma(j+1) - \Gamma(i+1). \tag{2.15}$$

Other non-elementary functions are also possible. In terms of the zero-th order Bessel function of the first kind (BFFK),

$$\varepsilon_{ij}(z) = \frac{J_0(i\lambda) - J_0(j\lambda)}{J_0(\lambda) - J_0(2\lambda)}; \ \lambda \in \mathbf{C} \tag{2.16}$$

a simple instance of which occurs for $\lambda = z_1 = 2.4048255576957727\cdots$, the first positive zero of the zero-th order BFFK, which yields the simpler expression

$$\varepsilon_{ij}(z_1) = \frac{J_0(z_1 j) - J_0(z_1 i)}{J_0(2z_1)}. \tag{2.17}$$

Since $z_1$ is irrational, (2.16, 17) are precision-dependent, like (2.13) which uses $\pi$, another irrational number. However, this is not problematic for most modern mathematical software.

Instead of orthogonal functions, orthogonal polynomials may be used, such as those of Hermite (*He*), of Laguerre (*L*), and of Jacobi, which are all feasible. Those of the latter include the Gegenbauer ($C^{(1)}$), the Chebyshev (*T*), and the Legendre (*P*) polynomials. They are all considered to be elementary functions, and in general,

$$\varepsilon_{ij}(z_1) = \frac{\Lambda(z_1 j) - \Lambda(z_1 i)}{\Lambda(2z_1)}; \ \Lambda \in \{He, L; C^{(1)}, P, T\}; \ i, j \in \{1, 2\}; \ z_1 \in \mathbf{R} \big| \Lambda(z_1) = 0 \tag{2.18}$$

for which $z_1$ may be the 1st positive zero of the orthogonal polynomial. For the 2nd order Laguerre polynomial for example,

$$L_2(z) = \left((z-2)^2 - 2\right)\big/2, \tag{2.19}$$

with $z_1 = 2 - \sqrt{2}$. Substituting (2.19) into (2.18) yields

$$\varepsilon_{ij}(z_1) = \frac{(z_1 j - 2)^2 - (z_1 i - 2)^2}{(2z_1 - 2)^2 - 2}. \tag{2.20}$$

The 2nd order polynomial (2.19) is only being used as an example. The first order Laguerre polynomial is equally effective, and yields an expression identical to (2.9). Eq. (2.10) may be adapted to higher dimensions without difficulty, which will be demonstrated in the following sections.

Another, different expression for (2.1) may be obtained using the following well-known integral definition of the Kronecker delta symbol:

$$\delta_{p-q} = \frac{1}{2\pi} \int_{-\pi}^{+\pi} e^{i(p-q)s} ds = \text{sinc}\left[\pi(p-q)\right], \ p, q \in \mathbf{Z}^+ \tag{2.21}$$

which vanishes everywhere, except where $p = q$. Applying this relation to (2.8) yields

$$\varepsilon_{ij} = \text{sinc}\left[\pi(j + 2i - 4)\right] - \text{sinc}\left[\pi(j - i + 1)\right]. \tag{2.22}$$

In this case, the numerator functions are dissimilar from those of the denominator, unlike (2.10). The numerator functions are the sine, whereas the denominator functions, the identity.

## 3. Three dimensional vector space

In three dimensions ($R^3$), the epsilon has 3 indices resulting $3^3 = 27$ distinct combinations, out of which 3! or 6 combinations have non-repetitive values, whereas the rest have at least 2 identical indices per combination:

$$\varepsilon_{ijk} = \begin{cases} +1 & \text{if } \{i,j,k\} \in \{\{1,2,3\},\{2,3,1\},\{3,1,2\}\}, \\ -1 & \text{if } \{i,j,k\} \in \{\{1,3,2\},\{2,1,3\},\{3,2,1\}\}, \\ 0 & \text{if } i=j, \text{ or } i=k, \text{ or } j=k. \end{cases} \quad (3.1)$$

When the set {1, 2, 3} is arranged in a circular configuration, the two non-zero cases of $\varepsilon_{ijk}$ are equivalent to clock-wise and counter-clockwise permutations [1].

The epsilon can also be expressed in terms of the signum function as

$$\varepsilon_{ijk} = \frac{j-i}{|j-i|} \frac{k-i}{|k-i|} \frac{k-j}{|k-j|} = \text{sgn}(j-i)\,\text{sgn}(k-i)\,\text{sgn}(k-j); \quad i,j,k \in \{1,2,3\} \quad (3.2)$$

which is a generalization of (2.2). When re-expressed in terms of the new index set $i_n$

$$\varepsilon_{i_1 i_2 i_3} = \frac{i_2-i_1}{|i_2-i_1|} \frac{i_3-i_1}{|i_3-i_1|} \frac{i_3-i_2}{|i_3-i_2|} = \prod_{3 \geq n \geq m \geq 1} \frac{i_n-i_m}{|i_n-i_m|} = \prod_{3 \geq n \geq m \geq 1} \text{sgn}(i_n - i_m); \quad i_1, i_2, i_3 \in \{1,2,3\}. \quad (3.3)$$

The bounds of the product symbol have to be correctly chosen in accordance with the inequality, prior to their use in the signum function. More explicitly however, the expression can be re-cast in terms of deltas as

$$\varepsilon_{i_1,i_2,i_3} = \prod_{m=1}^{3} \frac{i_{2+\delta(m-2)+\delta(m-3)} - i_{1+\delta(m-3)}}{|i_{2+\delta(m-2)+\delta(m-3)} - i_{1+\delta(m-3)}|}; \quad i_1, i_2, i_3 \in \{1,2,3\}. \quad (3.4)$$

Expressed differently, the subscript of the 1st index is 2 for $m = 1$, but is 3 otherwise. The subscript of the 2nd index is unity for all $m$ except at $m = 3$, for which it becomes 2. A conversion of the Kronecker delta symbol, derived in §**7** of [8], in now invoked,

$$\delta_{z-n} = [2\Gamma(z)\cos(n\pi) + z - 2](\Gamma(n) - 1) - (n\Gamma(z) - z)\cos(n\pi) + 1; \quad z, n \in \{1,2,3\} \quad (3.5)$$

which yields a delta for any $n \in \{1,2,3\}$ with $z$ constrained to {1, 2, 3}. Setting in turn $n = 2$ and 3 and substituting the 2 resultant expressions back into (3.4), a more explicit function than (3.3) is obtained for the first time, that eliminates the inequality in (3.3):

$$\varepsilon_{i_1,i_2,i_3} = \prod_{m=1}^{3} \frac{i_{m-\Gamma(m)+2} - i_{\Gamma(m)}}{|i_{m-\Gamma(m)+2} - i_{\Gamma(m)}|} = \prod_{m=1}^{3} \text{sgn}(i_{m-\Gamma(m)+2} - i_{\Gamma(m)}); \quad i_1, i_2, i_3 \in \{1,2,3\}. \quad (3.6)$$

If $i$ were considered to be the independent variable, $\varepsilon_{ijk}$ can be re-expressed as[1]

$$\varepsilon_{ijk} = \begin{cases} +1 & \text{if } j = (i+1)\delta_{i-1} + (i+1)\delta_{i-2} + (i-2)\delta_{i-3} \ \& \ k = (i+2)\delta_{i-1} + (i-1)\delta_{i-2} + (i-1)\delta_{i-3}, \\ -1 & \text{if } j = (i+2)\delta_{i-1} + (i-1)\delta_{i-2} + (i-1)\delta_{i-3} \ \& \ k = (i+1)\delta_{i-1} + (i+1)\delta_{i-2} + (i-2)\delta_{i-3}, \\ 0 & \text{if } i = j, \text{ or } i = k, \text{ or } j = k. \end{cases}$$

(3.7)

and after a simplification,

$$\varepsilon_{ijk} = \begin{cases} +1 & \text{if } j = 2\delta_{i-1} + 3\delta_{i-2} + \delta_{i-3} \ \& \ k = 3\delta_{i-1} + \delta_{i-2} + 2\delta_{i-3}, \\ -1 & \text{if } j = 3\delta_{i-1} + \delta_{i-2} + 2\delta_{i-3} \ \& \ k = 2\delta_{i-1} + 3\delta_{i-2} + \delta_{i-3}, \\ 0 & \text{if } i = j, \text{ or } i = k, \text{ or } j = k. \end{cases}$$

(3.8)

It can be recast as the equation

$$\varepsilon_{ijk} = \delta_{j-2\delta_{i-1}-3\delta_{i-2}-\delta_{i-3}} \delta_{k-3\delta_{i-1}-\delta_{i-2}-2\delta_{i-3}} - \delta_{j-3\delta_{i-1}-\delta_{i-2}-2\delta_{i-3}} \delta_{k-2\delta_{i-1}-3\delta_{i-2}-\delta_{i-3}}; \quad i,j,k \in \{1,2,3\}$$

(3.9)

The argument of a delta is now comprised of some linear combination of other deltas. It is possible to eliminate the deltas in the argument of each, using the observations

$$\delta_{i-1} = \frac{1}{2}(i-2)(i-3),$$  (3.10)

$$\delta_{i-2} = (1-i)(i-3),$$  (3.11)

$$\delta_{i-3} = \frac{1}{2}(i-1)(i-2),$$  (3.12)

and yielding the simpler expression

$$\varepsilon_{ijk} = \delta_{2j+3i^2-11i+4} \delta_{2k-3i^2+13i-16} - \delta_{2j-3i^2+13i-16} \delta_{2k+3i^2-11i+4}.$$  (3.13)

Both polynomial arguments have irrational and/or complex roots. However, a simpler expression can be obtained using the gamma function.

Using (3.5) again, with $n = 1, 2,$ and 3 and substituting all three resultant expressions into (3.9), yields in terms of the gamma function;

$$\varepsilon_{ijk} = \delta_{j+3\Gamma(i)-i-4} \delta_{k-3\Gamma(i)+2i-2} - \delta_{j-3\Gamma(i)+2i-2} \delta_{k+3\Gamma(i)-i-4} = \begin{vmatrix} \delta_{j+3\Gamma(i)-i-4} & \delta_{j-3\Gamma(i)+2i-2} \\ \delta_{k+3\Gamma(i)-i-4} & \delta_{k-3\Gamma(i)+2i-2} \end{vmatrix}.$$  (3.14)

---

[1] The variables $i$ and $j$ which are frequently used throughout this report, should not be confused with the imaginary number $i = (-1)^{1/2}$

It is known *a priori* that $\Gamma(1) = \Gamma(2) = 1$, whereas $\Gamma(3) = 2$ [9]. Furthermore, since the 2nd term is identical to the 1st term, but with the arguments of the delta symbols interchanged, (3.14) can also be expressed more compactly as

$$\varepsilon_{ijk} = \delta_{j'+3\Gamma(i)-i-4}\delta_{k'-3\Gamma(i)+2i-2}\left(\delta_{j'j}\delta_{k'k} - \delta_{j'k}\delta_{k'j}\right). \tag{3.15}$$

The minuend and the subtrahend are thus identical in form, but with the indices $j$ and $k$ interchanged. This expression is much more compact than one found from that due to Straub (1.7) [5], once it is specialized to $N = 3$.

The deltas may also be completely eliminated from (3.14). One method of doing so, is by using (3.5) again,

$$\varepsilon_{ijk} = \begin{matrix} \left[(2\Gamma(j)\cos(M\pi)+j-2)A-(M\Gamma(j)-j)\cos(M\pi)+1\right]\left[(2\Gamma(k)\cos(N\pi)+k-2)B-(N\Gamma(k)-k)\cos(N\pi)+1\right] \\ -\left[(2\Gamma(k)\cos(M\pi)+k-2)A-(M\Gamma(k)-k)\cos(M\pi)+1\right]\left[(2\Gamma(j)\cos(N\pi)+j-2)B-(N\Gamma(j)-j)\cos(N\pi)+1\right], \end{matrix} \tag{3.16}$$

where

$$A = \Gamma(M)-1,\ B = \Gamma(N)-1,\ M = -3\Gamma(i)+i+4,\ N = 3\Gamma(i)-2i+2, \tag{3.17}$$

and after some algebra, which is most easily reduced by simplifying the 1st term in (3.16), then subtracting it from its replica but with the indices $j$ and $k$ interchanged,

$$\varepsilon_{ijk} = \begin{aligned} &\left\{(2A-M)\left[(1+B(k-2))\Gamma(j)-(1+B(j-2))\Gamma(k)\right]+(j-k)(1-2B)\right\}\cos(\pi M) \\ &-\left\{(2B-N)\left[(1+A(k-2))\Gamma(j)-(1+A(j-2))\Gamma(k)\right]+(j-k)(1-2A)\right\}\cos(\pi N) \\ &+\left[(2(A-B)-M+N)(k\Gamma(j)-j\Gamma(k))\right]\cos(\pi M)\cos(\pi N)+(A-B)(j-k). \end{aligned} \tag{3.18}$$

Invoking (3.17),

$$A = \Gamma(M)-1 = \Gamma(-3\Gamma(i)+i+4)-1 = 1-2\Gamma(i)+i \tag{3.19}$$

$$B = \Gamma(N)-1 = \Gamma(3\Gamma(i)-2i+2)-1 = 1+\Gamma(i)-i \tag{3.20}$$

$$\cos(\pi M) = \cos\left[\pi(-3\Gamma(i)+i+4)\right] = 2\Gamma(i)-2i+1 \tag{3.21}$$

$$\cos(\pi N) = \cos\left[\pi(3\Gamma(i)-2i+2)\right] = 2\Gamma(i)-3 \tag{3.22}$$

which renders (3.18) as

$$\varepsilon_{ijk} = \begin{matrix} (j-k)\left[-12\Gamma^2(i)+(12i+7)\Gamma(i)-4(i^2+1)\right] \\ +\left(\Gamma(j)-\Gamma(k)\right)\left[12\Gamma^3(i)-(16i+2)\Gamma^2(i)+\left(12i^2-18i+12\right)\Gamma(i)-4i^3+12i^2-9i+2\right] \\ -\left(k\Gamma(j)-j\Gamma(k)\right)\left[6\Gamma^3(i)-(12i-7)\Gamma^2(i)+\left(10i^2-15i+2\right)\Gamma(i)-2i^3+i^2+8i-4\right] \end{matrix}$$

(3.23)

in which expressions in $i$ are decoupled from those in $j$ and in $k$. It can be seen that the epsilon vanishes for $j = k$, although other similar combinations are not so obvious. Since $i$ is constrained to the set $\{1, 2, 3\}$, the following simplifications are possible,

$$-12\Gamma^2(i)+(12i+7)\Gamma(i)-4(i^2+1)=-\Gamma(i) \qquad (3.24)$$

$$12\Gamma^3(i)-(16i+2)\Gamma^2(i)+\left(12i^2-18i+12\right)\Gamma(i)-4i^3+12i^2-9i+2=i \qquad (3.25)$$

$$6\Gamma^3(i)-(12i-7)\Gamma^2(i)+\left(10i^2-15i+2\right)\Gamma(i)-2i^3+i^2+8i-4=1 \qquad (3.26)$$

then in terms of gamma functions, the simplest expression is

$$\varepsilon_{ijk} = i\left(\Gamma(j)-\Gamma(k)\right)-\Gamma(i)(j-k)-\Gamma(j)k+j\Gamma(k). \qquad (3.27)$$

Using the identity,

$$\Gamma(z) = \frac{1}{2}\left(z^2-3z+4\right); \; z \in \{i,j,k\}; \; i,j,k \in \{1,2,3\}, \qquad (3.28)$$

then (3.27) simplifies to

$$\varepsilon_{ijk} = \frac{1}{2}(j-i)(k-i)(k-j) = \frac{1}{2}\prod_{m=1}^{3}\left(i_{2+m-\Gamma(m)}-i_{\Gamma(m)}\right); \; \{i,j,k\} \leftrightarrow \{i_1,i_2,i_3\} \qquad (3.29)$$

which is simpler than (3.6, 14). Eq. (3.29) is a special case of the more general relation

$$\varepsilon_{ijk}(\lambda) = \frac{G(i\lambda)-G(j\lambda)}{G(\lambda)-G(2\lambda)}\frac{G(i\lambda)-G(k\lambda)}{G(\lambda)-G(3\lambda)}\frac{G(j\lambda)-G(k\lambda)}{G(2\lambda)-G(3\lambda)}; \; i,j,k \in \{1,2,3\} \, \& \, \lambda \in \mathbf{C}. \qquad (3.30)$$

A compact version of this general expression may be obtained using the gamma function, after assigning the original indices $\{i, j, k\}$ bijectively to the new indices $\{i_1, i_2, i_3\}$,

$$\varepsilon_{i_1 i_2 i_3}(\lambda) = \prod_{n=1}^{3}\frac{G\left(i_{\Gamma(n)}\lambda\right)-G\left(i_{2+n-\Gamma(n)}\lambda\right)}{G\left(\Gamma(n)\lambda\right)-G\left((2+n-\Gamma(n))\lambda\right)}; \; i_1,i_2,i_3 \in \{1,2,3\} \, \& \, \lambda \in \mathbf{C}. \qquad (3.31)$$

If the general function is set to the identity or $G(\xi) = \xi$, for instance, (3.30) yields (3.29)

$$\varepsilon_{ijk} = \frac{i\lambda - j\lambda}{\lambda - 2\lambda} \frac{i\lambda - k\lambda}{\lambda - 3\lambda} \frac{j\lambda - k\lambda}{2\lambda - 3\lambda} = \frac{(i-j)(i-k)(j-k)}{-1 \times -2 \times -1} = \frac{1}{2}(j-i)(k-i)(k-j). \quad (3.32)$$

For a trigonometric function, with $G(\xi) = \cos(\xi)$ in (3.30),

$$\varepsilon_{ijk}(\lambda) = \frac{\cos(i\lambda) - \cos(j\lambda)}{\cos(\lambda) - \cos(2\lambda)} \frac{\cos(i\lambda) - \cos(k\lambda)}{\cos(\lambda) - \cos(3\lambda)} \frac{\cos(j\lambda) - \cos(k\lambda)}{\cos(2\lambda) - \cos(3\lambda)} \quad (3.33)$$

a simple case of which is that for $\lambda = \pi/4$,

$$\varepsilon_{ijk}\left(\frac{\pi}{4}\right) = \sqrt{2}\left(\cos\left(\frac{i\pi}{4}\right) - \cos\left(\frac{j\pi}{4}\right)\right)\left(\cos\left(\frac{k\pi}{4}\right) - \cos\left(\frac{i\pi}{4}\right)\right)\left(\cos\left(\frac{k\pi}{4}\right) - \cos\left(\frac{j\pi}{4}\right)\right) \quad (3.34)$$

which can also have other forms using well-known trigonometric identities,

$$\varepsilon_{ijk}\left(\frac{\pi}{4}\right) = 2\sin\left(\frac{(j-i)\pi}{4}\right)\sin\left(\frac{(k-i)\pi}{4}\right)\sin\left(\frac{(k-j)\pi}{4}\right) = 2\prod_{n=1}^{3}\sin\left(\frac{(i_{2+n-\Gamma(n)} - i_{\Gamma(n)})\pi}{4}\right). \quad (3.35)$$

Unlike the LHS of the equation, the RHS is no longer in terms of elementary functions due to the use of the gamma function.

In terms of a non-elementary function such as the gamma function, the simplest expression is obtained for $\lambda = 1$ but with $G(\xi) = \Gamma(\xi + 1)$ instead of $G(\xi) = \Gamma(\xi)$ in (3.30),

$$\varepsilon_{ijk} = \frac{1}{20}\left(\Gamma(j+1) - \Gamma(i+1)\right)\left(\Gamma(k+1) - \Gamma(i+1)\right)\left(\Gamma(k+1) - \Gamma(j+1)\right), \quad (3.36)$$

and can be simplified to the compact expression after the assignment $\{i, j, k\} \leftrightarrow \{i_1, i_2, i_3\}$,

$$\varepsilon_{ijk} = \frac{1}{20}\prod_{n=1}^{3}\left[\Gamma(i_{2+n-\Gamma(n)} + 1) - \Gamma(i_{\Gamma(n)} + 1)\right]. \quad (3.37)$$

In terms of the zero-th order Bessel function of the first kind (BFFK), $G(\xi) = J_0(\xi)$,

$$\varepsilon_{ijk}(\lambda) = \frac{J_0(\lambda i) - J_0(\lambda j)}{J_0(\lambda) - J_0(2\lambda)} \frac{J_0(\lambda i) - J_0(\lambda k)}{J_0(\lambda) - J_0(3\lambda)} \frac{J_0(\lambda j) - J_0(\lambda k)}{J_0(2\lambda) - J_0(3\lambda)} \quad (3.38)$$

but with $\lambda = z_1 = 2.4048255576957727\cdots$, the 1st zero of the BFFK, (3.38) simplifies to

$$\varepsilon_{ijk} = \frac{\left(J_0(z_1 i) - J_0(z_1 j)\right)\left(J_0(z_1 i) - J_0(z_1 k)\right)\left(J_0(z_1 j) - J_0(z_1 k)\right)}{J_0(2z_1) J_0(3z_1)\left(J_0(2z_1) - J_0(3z_1)\right)} \quad . \tag{3.39}$$

The above expression is not limited to the zero-th order BFFK.

As in the previous section, this section's epsilon is also expressible in terms of other orthogonal functions or polynomials $\Lambda$. In general, it can be deduced that

$$\varepsilon_{ijk}(z_1) = \frac{\left(\Lambda(z_1 j) - \Lambda(z_1 i)\right)\left(\Lambda(z_1 k) - \Lambda(z_1 i)\right)\left(\Lambda(z_1 k) - \Lambda(z_1 j)\right)}{\Lambda(2z_1)\Lambda(3z_1)\left(\Lambda(2z_1) - \Lambda(3z_1)\right)} ; \tag{3.40}$$

$$\Lambda \in \{He, L; C^{(1)}, P, T\}; \ i, j, k \in \{1, 2, 3\}; \ z_1 \in \mathbb{R} \,|\, \Lambda(z_1) = 0$$

which, like (3.39), is in the same form as (3.30). The appropriate trigonometric and Bessel functions, such as the ones used in this section, and which are orthogonal functions, also satisfy the above expression. The denominator is simplified as in the above expression only if $z_1$ is selected to be a zero of $\Lambda(z)$.

As in the previous section, using the integral identity (2.21) of the Kronecker delta symbol, it is possible to obtain other expressions not of the form of (3.30, 31). Applying (2.21) to (3.13) yields

$$\varepsilon_{ijk} = \begin{array}{l} \operatorname{sinc}\left[\pi(2j + 3i^2 - 11i + 4)\right] \operatorname{sinc}\left[\pi(2k - (3i^2 - 13i + 16))\right] \\ - \operatorname{sinc}\left[\pi(2k + 3i^2 - 11i + 4)\right] \operatorname{sinc}\left[\pi(2j - (3i^2 - 13i + 16))\right] \end{array} \tag{3.41}$$

which is wholly in terms of elementary functions. It can be expressed as a 2 x 2 determinant. It can also be simplified to the compact expression

$$\varepsilon_{ijk} = \sum_{n=2}^{3} (-1)^n \operatorname{sinc}\left[\pi(2i_n + 3i_1^2 - 11i_1 + 4)\right] \operatorname{sinc}\left[\pi(2i_{5-n} - (3i_1^2 - 13i_1 + 16))\right]. \tag{3.42}$$

Another expression may be obtained by applying (2.21) to (3.14) instead, which yields,

$$\varepsilon_{ijk} = \begin{array}{l} \operatorname{sinc}\left[\pi(j + 3\Gamma(i) - i - 4)\right] \operatorname{sinc}\left[\pi(k - 3\Gamma(i) + 2i - 2)\right] \\ - \operatorname{sinc}\left[\pi(k + 3\Gamma(i) - i - 4)\right] \operatorname{sinc}\left[\pi(j - 3\Gamma(i) + 2i - 2)\right] \end{array} \tag{3.43}$$

or

$$\varepsilon_{ijk} = \sum_{n=2}^{3} (-1)^n \operatorname{sinc}\left[\pi(i_n + 3\Gamma(i) - i - 4)\right] \operatorname{sinc}\left[\pi(i_{5-n} - 3\Gamma(i) + 2i - 2)\right] \tag{3.44}$$

which is simpler and more compact than (3.42, 43). However, using (2.21) yields increasingly cumbersome expressions for higher dimensions, and is not used beyond this section. By contrast, (3.30, 31) and (3.40) can be more easily extended to higher dimensions, while retaining most of their compactness.

## 4. Applications

The following well-known identity in $R^2$,

$$\varepsilon_{ij}\varepsilon_{mn} = \delta_{im}\delta_{jn} - \delta_{in}\delta_{jm} \; ; \; i,j,m,n \in \{1,2\} \tag{4.1}$$

can be re-expressed in terms of elementary functions using (2.11),

$$\varepsilon_{ij}\varepsilon_{mn} = (j-i)(n-m), \tag{4.2}$$

which is a much simpler, and a more memorable expression than (4.1).

In $R^3$, the following expression is the result of a 3 x 3 determinant:

$$\varepsilon_{ijk}\varepsilon_{lmn} = \delta_{il}\left(\delta_{jm}\delta_{kn} - \delta_{jn}\delta_{km}\right) - \delta_{im}\left(\delta_{jl}\delta_{kn} - \delta_{jn}\delta_{kl}\right) + \delta_{in}\left(\delta_{jl}\delta_{km} - \delta_{jm}\delta_{kl}\right);$$
$$i,j,k,l,m,n \in \{1,2,3\}. \tag{4.3}$$

Invoking (3.29), a more compact, and a simpler analytical expression results, entirely in terms of elementary functions,

$$\varepsilon_{ijk}\varepsilon_{lmn} = \frac{1}{4}(j-i)(k-i)(k-j)(m-l)(n-l)(n-m) \tag{4.4}$$

and which can also be expressed in terms of the gamma function as

$$\varepsilon_{ijk}\varepsilon_{lmn} = \frac{1}{4}\prod_{p=1}^{3}\left(i_{p-\Gamma(p)+2} - i_{\Gamma(p)}\right)\left(i_{p-\Gamma(p)+5} - i_{\Gamma(p)+3}\right); \; \{i,j,k,l,m,n\} \leftrightarrow \{i_1,i_2,i_3,i_4,i_5,i_6\}. \tag{4.5}$$

It is more compact than either (4.3) or (4.4), although it is no longer in terms of elementary functions due to the use of the gamma function. However, it is known *a priori* that $\Gamma(1) = \Gamma(2) = 1$, and that $\Gamma(3) = 2$ [9].

The determinant of a 3 x 3 matrix $\mathbf{X}$, which is widely used in vector analysis, is known to have the following expression in terms of the three-dimensional epsilon [1],

$$|\mathbf{X}| = \sum_{i=1}^{3}\sum_{j=1}^{3}\sum_{k=1}^{3} \varepsilon_{ijk} x_{1,i} x_{2,j} x_{3,k} \tag{4.6}$$

which is somewhat difficult to expand without consulting or memorizing (3.1), or beginning with {1, 2, 3} and interchanging adjacent values until the desired permutation is attained. However, a more analytical expression can be found by using (3.29) again,

$$|\mathbf{X}| = \frac{1}{2}\sum_{i=1}^{3}\sum_{j=1}^{3}\sum_{k=1}^{3}(j-i)(k-i)(k-j) x_{1,i} x_{2,j} x_{3,k} = \frac{1}{2}\sum_{i_1=1}^{3}\sum_{i_2=1}^{3}\sum_{i_3=1}^{3}\prod_{m=1}^{3}\left(i_{m-\Gamma(m)+2} - i_{\Gamma(m)}\right) x_{m,i_m}. \tag{4.7}$$

Although some compactness is sacrificed, the expression is more amenable to an expeditious substitution of values using the summations, than is (4.6).

## 5. Summary and Conclusions

The Levi-Civita symbol (or epsilon) in $R^2$, is known to be expressible as [1]

$$\varepsilon_{ij} = \text{sgn}(j-i) = \begin{cases} +1 & \text{if } \{i,j\} = \{1,2\}, \\ -1 & \text{if } \{i,j\} = \{2,1\}, \\ 0 & \text{if } \{i,j\} \in \{\{1,1\},\{2,2\}\} \end{cases} \tag{5.1}$$

but was found to have the simpler alternative

$$\varepsilon_{ij} = j - i \tag{5.2}$$

which although easily deduced by inspection from (5.1), was derived beginning with (5.1), in §**2**.

In $R^3$, it is known that [1]

$$\varepsilon_{ijk} = \prod_{3 \geq n \geq m \geq 1} \text{sgn}(i_n - i_m) = \begin{cases} +1 & \text{if } \{i,j,k\} \in \{\{1,2,3\},\{2,3,1\},\{3,1,2\}\}, \\ -1 & \text{if } \{i,j,k\} \in \{\{1,3,2\},\{2,1,3\},\{3,2,1\}\}, \\ 0 & \text{if } i = j, \text{ or } i = k, \text{ or } j = k. \end{cases} \tag{5.3}$$

As previously remarked, (5.1) and (5.3) have more in common with a look-up table than to actual equations. In §**3** and for the first time, a more explicit expression for the LHS of (5.3) that eliminates the inequality on the product symbol, has been found to be

$$\varepsilon_{i_1,i_2,i_3} = \prod_{m=1}^{3} \frac{i_{m-\Gamma(m)+2} - i_{\Gamma(m)}}{|i_{m-\Gamma(m)+2} - i_{\Gamma(m)}|} = \prod_{m=1}^{3} \text{sgn}(i_{m+2-\Gamma(m)} - i_{\Gamma(m)}); \quad i_1, i_2, i_3 \in \{1,2,3\} \tag{5.4}$$

and makes use of the gamma function, which is known to be continuous for positive, integral arguments [9]. The use of the gamma function should present no difficulty, as it is known that $\Gamma(1) = \Gamma(2) = 1$, while $\Gamma(3) = 2$ [9]. In terms of Kronecker delta symbols (or deltas), the following simple expression has also been found for the first time:

$$\varepsilon_{ijk} = \delta_{j+3\Gamma(i)-i-4}\delta_{k-3\Gamma(i)+2i-2} - \delta_{j-3\Gamma(i)+2i-2}\delta_{k+3\Gamma(i)-i-4}; \quad i,j,k \in \{1,2,3\}. \tag{5.5}$$

However, the signum function, the gamma function, and the delta symbol, are not elementary functions according to the Risch Algorithm [2-4]. The following relation was also derived in §**3**, beginning with (5.3),

$$\varepsilon_{ijk} = \frac{1}{2}(j-i)(k-i)(k-j) = \frac{1}{2}\prod_{n=1}^{3}(i_{n+2-\Gamma(n)} - i_{\Gamma(n)}) = \varepsilon_{i_1 i_2 i_3} \tag{5.6}$$

whose LHS is now in terms of elementary functions, and which can also be found by a mere inspection of the RHS of (5.3). The RHS expression is however reported here for the first time. The three-dimensional epsilon is the only epsilon for which the number of product product terms is identical to the number of indices.

In $R^4$, the description of the epsilon becomes significantly larger. The epsilon now has four indices resulting in $4^4 = 256$ distinct combinations, out of which 4! or 24 combinations have non-repetitive indices, whereas the rest have at least 2 identical indices for each combination,

$$\varepsilon_{ijkl} = \begin{cases} +1 & \text{if } \{i,j,k,l\} \in \begin{Bmatrix} \{1,2,3,4\},\{1,3,4,2\},\{1,4,2,3\};\{2,1,4,3\},\{2,3,1,4\},\{2,4,3,1\}; \\ \{3,1,2,4\},\{3,2,4,1\},\{3,4,1,2\};\{4,1,3,2\},\{4,2,1,3\},\{4,3,2,1\} \end{Bmatrix}, \\ -1 & \text{if } \{i,j,k,l\} \in \begin{Bmatrix} \{1,3,2,4\},\{1,4,3,2\},\{1,2,4,3\};\{2,4,1,3\},\{2,1,3,4\},\{2,3,4,1\}; \\ \{3,2,1,4\},\{3,4,2,1\},\{3,1,4,2\};\{4,3,1,2\},\{4,1,2,3\},\{4,2,3,1\} \end{Bmatrix}, \\ 0 & \text{otherwise.} \end{cases}$$

(5.7)

However, an analytical and a more succinct expression can be deduced to be

$$\varepsilon_{ijkl} = \frac{1}{12}(j-i)(k-i)(l-i)(k-j)(l-j)(l-k); \quad i,j,k,l \in \{1,2,3,4\}.$$  (5.8)

It is the first epsilon for which the number of product terms (6) exceeds the number of its indices (4). It can be re-cast more compactly as

$$\varepsilon_{i_1 i_2 i_3 i_4} = \frac{1}{12} \prod_{n=1}^{3} (i_{n+1} - i_1)(i_{n+3-\Gamma(n)} - i_{\Gamma(n)+1}); \quad i_1, i_2, i_3, i_4 \in \{1,2,3,4\}.$$  (5.9)

For $N \geq 5$, it becomes increasingly impractical to express the epsilon as a collection of explicit conditions, like (5.7). Furthermore, the contracting procedure used for (5.8, 9) yields increasingly larger products, prompting the search for another, more efficient approach. One such alternative is the following, rational expression

$$\varepsilon_{i_1 i_2 i_3 \cdots i_N} = \prod_{m=1}^{N-1} \prod_{n=1}^{N-m} \frac{i_{N+1-m} - i_n}{N+1-m-n}; \quad i_1, i_2, i_3, \cdots, i_N \in \{1,2,3,\cdots,N\}, \; N > 1.$$  (5.10)

In this expression, which is entirely in terms of elementary functions, the inner or nested $n$-product is not known, unless the index $m$ is specified by the outer $m$-product. For $N = 5$ for instance,

$$\varepsilon_{i_1 i_2 i_3 i_4 i_5} = \prod_{n=1}^{4} \frac{i_5 - i_n}{5-n} \prod_{n=1}^{3} \frac{i_4 - i_n}{4-n} \prod_{n=1}^{2} \frac{i_3 - i_n}{3-n} \prod_{n=1}^{1} \frac{i_2 - i_n}{2-n}$$

$$= \frac{1}{288}(i_2-i_1)(i_3-i_1)(i_4-i_1)(i_5-i_1)(i_3-i_2)(i_4-i_2)(i_5-i_2)(i_4-i_3)(i_5-i_3)(i_5-i_4).$$

(5.11)

It should be obvious that if any pair of indices are identical in a parenthesized difference, the epsilon vanishes.

Lastly, it was found in sections **2** and **3** that an epsilon can also be re-expressed using different, and properly normalized, non-elementary (special) functions. Since (5.10) encapsulates the cases for which $N = 2$ and $3$, which were the ones examined in sections **2** and **3**, it is not difficult to conclude that this re-expression should extend to (5.10). As a generalized, discrete function, the epsilon may be re-expressed as follows

$$\varepsilon_{i_1 i_2 i_3 \cdots i_N}(\lambda) = \prod_{m=1}^{N-1} \prod_{n=1}^{N-m} \frac{G(i_{N+1-m}\lambda) - G(i_n \lambda)}{G((N+1-m)\lambda) - G(n\lambda)} \; ; i_1, i_2, i_3, \cdots, i_N \in \{1, 2, 3, \cdots, N\}; N > 1; \lambda \in \mathbf{C}$$
(5.12)

with $\lambda$ being generally complex. Examples of (5.12) *respectively* include, but are not limited to, the cosine function, the zero-th order Bessel function of the first kind (BFFK), the gamma function, and (non-)orthogonal polynomials $\Lambda$,

$$\varepsilon_{i_1 i_2 i_3 \cdots i_N}(\lambda) = \prod_{m=1}^{N-1} \prod_{n=1}^{N-m} \frac{\cos(i_{N+1-m}\lambda) - \cos(i_n \lambda)}{\cos((N+1-m)\lambda) - \cos(n\lambda)},$$
(5.13)

$$\varepsilon_{i_1 i_2 i_3 \cdots i_N}(\lambda) = \prod_{m=1}^{N-1} \prod_{n=1}^{N-m} \frac{J_0(i_{N+1-m}\lambda) - J_0(i_n \lambda)}{J_0((N+1-m)\lambda) - J_0(n\lambda)},$$
(5.14)

$$\varepsilon_{i_1 i_2 i_3 \cdots i_N}(1) = \prod_{m=1}^{N-1} \prod_{n=1}^{N-m} \frac{\Gamma(i_{N+1-m}+1) - \Gamma(i_n+1)}{\Gamma(N+2-m) - \Gamma(n+1)},$$
(5.15)

$$\varepsilon_{i_1 i_2 i_3 \cdots i_N}(\lambda) = \prod_{m=1}^{N-1} \prod_{n=1}^{N-m} \frac{\Lambda(i_{N+1-m}\lambda) - \Lambda(i_n \lambda)}{\Lambda((N+1-m)\lambda) - \Lambda(n\lambda)}.$$
(5.16)

Eq. (5.13) is not restricted to the cosine function, as other trigonometric functions may also be used. Likewise, (5.14) is not restricted to the zero-th order BFFK, and other Bessel functions can be equally viable. Instead of orthogonal functions such as trigonometric and Bessel functions, (non-)orthogonal polynomials $\Lambda$ may also be used. They include those of Hermite (*He*), of Laguerre (*L*), and of Jacobi. Those of Jacobi include the Gegenbauer ($C^{(1)}$), the Chebyshev (*T*), and the Legendre (*P*) polynomials. Unlike Bessel and gamma functions, they are all considered to be elementary functions, like the cosine function. Other, polynomial candidates for $\Lambda$ include the Bell [10, 11] and the Bernoulli polynomials [13, 14], which are also elementary. The parameter $\lambda$ may in general be complex, with one exception. For (5.15), $\lambda$ is required to be unity in order for (5.15) to replicate the behavior of the epsilon in $\mathbf{R}^N$. The above equations represent just a few of many possibilities, but include both elementary and non-elementary functions. Eqs. (5.14) and (5.15) represent non-elementary functions, according to the Risch Algorithm [2-4]. The simplest expression in terms of elementary functions is (5.10), and corresponds to the identity ($G(\xi) = \xi$) function, and without restrictions on $\lambda$.

Other expressions are possible in terms of infinite series [15], or in terms of integral representations of the Kronecker delta symbol [16, 17]. They arise as a consequence of the orthogonality of basis functions. Although they may merit an investigation, they are unlikely to yield similarly compact expressions, and were therefore not considered in this report.

Three appendices follow the **References** section. In **Appendix A**, an implementation of (5.11) is encoded in Matlab. In **Appendix B**, a Matlab program is used to demonstrate that the epsilon can be treated as a generalized, discrete function, using (5.11, 13-15), all of which are confirmed to produce identical results for $N = 5$. **Appendix C** replicates the program of **Appendix A** using a vectorization of for-loops.

## Acknowledgements

Software such as Matlab®, Maple®, and Maxima were found useful in verifying some of the derived expressions. www.WolframAlpha.com was also helpful in this regard, and is widely accessible. Lastly, this document was generated in Microsoft Word, and converted prior to the upload to arXiv.org using www.freepdfconvert.com

# Appendix A

The following Matlab code generates the epsilon $\varepsilon_{i_1 i_2 i_3 i_4 i_5}$ (5.11) in $R^5$, but only displays the combinations $\{i_1, i_2, i_3, i_4, i_5\}$ for which the epsilon is not zero. The sole if-statement may be disabled to also display all the zero results. Nested for-loops are not considered to be efficient programming in general, but are used in this case to maintain simplicity and accessibility, although not compactness. A more compact version of this code using vectorized for-loops is found in **Appendix C**, however. The code below may be copied and pasted into the command window of Matlab for a quick execution

```matlab
clc, clear all

s1='                   ';s2='                ';s3='           ';
disp([s1,'Only the Non-Zero Cases of the Epsilon Will Be Displayed',sprintf('\n')])
disp(['the Enter-key must be used to advance the computation for each non-zero {i1, i2, i3, i4, i5} set:',sprintf('\n')])

N=5; % assuming the Levi-Civita Symbol (epsilon) in R^5 Space

M=0; % initialize epsilon counter
% compute all combinations of the epsilon
for i1=1:1:N
   for i2=1:1:N
      for i3=1:1:N
         for i4=1:1:N
            for i5=1:1:N,
    epsilon=(i2-i1)*(i3-i1)*(i4-i1)*(i5-i1)*(i3-i2)*(i4-i2)*(i5-i2) *(i4-i3)*(i5-i3)*(i5-i4)./ ...
            [(2- 1)*(3 -1)*(4 - 1)*( 5- 1)*(3 - 2)*(4 -2)*( 5- 2)*( 4- 3)*(5 - 3)*(5 - 4)];
              if epsilon ~= 0 % only display the non-zero values of epsilon:
disp(strcat([s1,s2,'{i1, i2, i3, i4, i5} = {'],num2str(i1),',',num2str(i2),',',num2str(i3),',',num2str(i4),',',num2str(i5),'}',[', epsilon =' ' ' num2str(epsilon)],';'));
              M=M+1; % number of unique non-zero epsilons
              pause
              end
            end
         end
      end
   end
end

disp(strcat([sprintf('\n'),s3,s3,s3,'Total number of non-zero epsilons = '], num2str(M),' or 5!'))
```

# Appendix B

The following Matlab code generates the epsilon $\varepsilon_{i_1 i_2 i_3 i_4 i_5}$ in $R^5$, and also displays all the combinations $\{i_1, i_2, i_3, i_4, i_5\}$ for which the epsilon is not zero. It also demonstrates that the epsilon can be treated as a generalized, discrete function, as explained for (5.13-15) in § **5**

```matlab
s1='            ';s2='          ';s3='        ';
disp([s1,'Only the Non-Zero Cases of the Epsilon Will Be Displayed', sprintf('\n')])
disp(['the Enter-key must be used to advance the computation for each non-zero {i1, i2, i3, i4, i5} set',sprintf('\n')])
METHOD=input([s1,'select method (e.g. Bessel, or cosine, or gamma, or the identity (default)) : '],'s');
K=0; N=5; % assuming the Levi-Civita Symbol (epsilon) in R^5 Space
lambda=randn+sqrt(-1)*randn; % random complex parameter
% initialize inner n-product and outer m-product variables:
iProduct=1; oProduct=1;
% initialize each of the N, x-indices, to {1,2,3,...,N}:
x=meshgrid(1:N);

for i1=1:1:N
 for i2=1:1:N
  for i3=1:1:N
   for i4=1:1:N
    for i5=1:1:N
       p=[i1 i2 i3 i4 i5];
       oProduct=1;
       for m=1:1:N-1
         iProduct=1;
         for n=1:1:N-m
           switch lower(METHOD)
             case 'bessel'
                iProduct=iProduct*[bessel(0,lambda*x(N+1-m, p(N+1-m)))-bessel(0, lambda*x(n,p(n)))]./ ...
                          [bessel(0, lambda* (N+1-m))-bessel(0, lambda*n)]; % Eq. (5.14)
             case 'cosine'
                iProduct=iProduct*[cos(lambda*x(N+1-m, p(N+1-m)))-cos(lambda*x(n, p(n)))]./ ...
                          [cos(lambda* (N+1-m))-cos(lambda*n)]; % Eq. (5.13)
             case 'gamma'
                iProduct=iProduct*[gamma(x(N+1-m,p(N+1-m))+1)-gamma(x(n,p(n))+1)]./ ...
                          [gamma((N+1-m)+1)-gamma(n+1)]; % Eq. (5.15)
             otherwise
                iProduct=iProduct*[x(N+1-m, p(N+1-m))-x(n, p(n))]./[(N+1-m)-n]; % Eq. (5.10) or (5.11)
           end % switch
         end % n
         oProduct=oProduct*iProduct; % Eq. (5.12)
       end % m
       epsilon=round(oProduct); % eliminates precision-induced imaginary part of order of 10^(-16)
        if epsilon ~=0 % only display the non-zero epsilons:
          disp([s1,s2,'{i1, i2, i3, i4, i5} = {',
              num2str(i1),',',num2str(i2),',',num2str(i3),',',num2str(i4),',',num2str(i5),'}', ...
              ', epsilon =' ' ' num2str(epsilon),';']);
         pause
         K=K+1;
        end % if
    end
   end
  end
 end
end

disp([sprintf('\n'),s3,s3,s3,'Total number of non-zero epsilons = ', num2str(K),' or 5!'])
```

# Appendix C

The following Matlab code replicates the results of **Appendix A**, but with a vectorization of for-loops instead, resulting in a much more compact code, although slightly more esoteric. Furthermore, the execution time is reduced by a factor of $\approx 10$. The computer hosting the Matlab software was the same one used for [8]. This code, like those in the preceding appendices, was generated using Matlab R2009, with which is available subroutines such as ngrid() and ind2sub() that are used in the code. The results were cross-checked at www.Wolframalpha.com using the Mathematica function Signature[{i1, i2, i3, i4, i5}], which is an implementation of the Levi-Civita Symbol according to https://mathworld.wolfram.com/PermutationSymbol.html

```
disp(sprintf('\n\n'))
[i1,i2,i3,i4,i5] = ndgrid(1:5,1:5,1:5,1:5,1:5);
epsilon=(1/288)*(i2-i1).*(i3-i1).*(i4-i1).*(i5-i1).*(i3-i2).*(i4-i2).*(i5-i2).*(i4-i3) ...
    .*(i5-i3).*(i5-i4); % Eq. (5.11)
[i_,j_,k_,l_,m_] = ind2sub(size(epsilon),find(epsilon == -1));% only the -1's
disp(['{i1,i2,i3,i4,i5} combinations that yield epsilon = -1:',sprintf('\n')])
disp([i_ j_ k_ l_ m_]), disp(sprintf('\n'))
[i,j,k,l,m] = ind2sub(size(epsilon),find(epsilon == 1));% only the +1's
disp(['{i1,i2,i3,i4,i5} combinations that yield epsilon = +1:',sprintf('\n')])
disp([i j k l m])
```

20